\documentclass{amsart}

\usepackage{amsfonts,amssymb,verbatim,amsmath,amsthm,latexsym,textcomp,amscd}
\usepackage{latexsym,amsfonts,amssymb,epsfig,verbatim}
\usepackage{amsmath,amsthm,amssymb,latexsym,graphics,textcomp}
\usepackage{graphicx}
\usepackage{color}
\usepackage{url}

\input xy
\xyoption{all}

\begin{document}

\newtheorem{theorem}{Theorem}[section]
\newtheorem{prop}[theorem]{Proposition}
\newtheorem{lemma}[theorem]{Lemma}
\newtheorem{cor}[theorem]{Corollary}
\newtheorem{definition}[theorem]{Definition}
\newtheorem{conj}[theorem]{Conjecture}
\newtheorem{rmk}[theorem]{Remark}
\newtheorem{claim}[theorem]{Claim}
\newtheorem{defth}[theorem]{Definition-Theorem}

\newcommand{\boundary}{\partial}
\newcommand{\C}{{\mathbb C}}
\newcommand{\integers}{{\mathbb Z}}
\newcommand{\natls}{{\mathbb N}}
\newcommand{\ratls}{{\mathbb Q}}
\newcommand{\reals}{{\mathbb R}}
\newcommand{\proj}{{\mathbb P}}
\newcommand{\lhp}{{\mathbb L}}
\newcommand{\tube}{{\mathbb T}}
\newcommand{\cusp}{{\mathbb P}}
\newcommand\AAA{{\mathcal A}}
\newcommand\BB{{\mathcal B}}
\newcommand\CC{{\mathcal C}}
\newcommand\DD{{\mathcal D}}
\newcommand\EE{{\mathcal E}}
\newcommand\FF{{\mathcal F}}
\newcommand\GG{{\mathcal G}}
\newcommand\HH{{\mathcal H}}
\newcommand\II{{\mathcal I}}
\newcommand\JJ{{\mathcal J}}
\newcommand\KK{{\mathcal K}}
\newcommand\LL{{\mathcal L}}
\newcommand\MM{{\mathcal M}}
\newcommand\NN{{\mathcal N}}
\newcommand\OO{{\mathcal O}}
\newcommand\PP{{\mathcal P}}
\newcommand\QQ{{\mathcal Q}}
\newcommand\RR{{\mathcal R}}
\newcommand\SSS{{\mathcal S}}
\newcommand\TT{{\mathcal T}}
\newcommand\UU{{\mathcal U}}
\newcommand\VV{{\mathcal V}}
\newcommand\WW{{\mathcal W}}
\newcommand\XX{{\mathcal X}}
\newcommand\YY{{\mathcal Y}}
\newcommand\ZZ{{\mathcal Z}}
\newcommand\CH{{\CC\HH}}
\newcommand\PEY{{\PP\EE\YY}}
\newcommand\MF{{\MM\FF}}
\newcommand\RCT{{{\mathcal R}_{CT}}}
\newcommand\PMF{{\PP\kern-2pt\MM\FF}}
\newcommand\FL{{\FF\LL}}
\newcommand\PML{{\PP\kern-2pt\MM\LL}}
\newcommand\GL{{\GG\LL}}
\newcommand\Pol{{\mathcal P}}
\newcommand\half{{\textstyle{\frac12}}}
\newcommand\Half{{\frac12}}
\newcommand\Mod{\operatorname{Mod}}
\newcommand\Area{\operatorname{Area}}
\newcommand\ep{\epsilon}
\newcommand\hhat{\widehat}
\newcommand\Proj{{\mathbf P}}
\newcommand\U{{\mathbf U}}
 \newcommand\Hyp{{\mathbf H}}
\newcommand\D{{\mathbf D}}
\newcommand\Z{{\mathbb Z}}
\newcommand\R{{\mathbb R}}
\newcommand\Q{{\mathbb Q}}
\newcommand\E{{\mathbb E}}
\newcommand\til{\widetilde}
\newcommand\length{\operatorname{length}}
\newcommand\tr{\operatorname{tr}}
\newcommand\gesim{\succ}
\newcommand\lesim{\prec}
\newcommand\simle{\lesim}
\newcommand\simge{\gesim}
\newcommand{\simmult}{\asymp}
\newcommand{\simadd}{\mathrel{\overset{\text{\tiny $+$}}{\sim}}}
\newcommand{\ssm}{\setminus}
\newcommand{\diam}{\operatorname{diam}}
\newcommand{\pair}[1]{\langle #1\rangle}
\newcommand{\T}{{\mathbf T}}
\newcommand{\inj}{\operatorname{inj}}
\newcommand{\pleat}{\operatorname{\mathbf{pleat}}}
\newcommand{\short}{\operatorname{\mathbf{short}}}
\newcommand{\vertices}{\operatorname{vert}}
\newcommand{\collar}{\operatorname{\mathbf{collar}}}
\newcommand{\bcollar}{\operatorname{\overline{\mathbf{collar}}}}
\newcommand{\I}{{\mathbf I}}
\newcommand{\tprec}{\prec_t}
\newcommand{\fprec}{\prec_f}
\newcommand{\bprec}{\prec_b}
\newcommand{\pprec}{\prec_p}
\newcommand{\ppreceq}{\preceq_p}
\newcommand{\sprec}{\prec_s}
\newcommand{\cpreceq}{\preceq_c}
\newcommand{\cprec}{\prec_c}
\newcommand{\topprec}{\prec_{\rm top}}
\newcommand{\Topprec}{\prec_{\rm TOP}}
\newcommand{\fsub}{\mathrel{\scriptstyle\searrow}}
\newcommand{\bsub}{\mathrel{\scriptstyle\swarrow}}
\newcommand{\fsubd}{\mathrel{{\scriptstyle\searrow}\kern-1ex^d\kern0.5ex}}
\newcommand{\bsubd}{\mathrel{{\scriptstyle\swarrow}\kern-1.6ex^d\kern0.8ex}}
\newcommand{\fsubeq}{\mathrel{\raise-.7ex\hbox{$\overset{\searrow}{=}$}}}
\newcommand{\bsubeq}{\mathrel{\raise-.7ex\hbox{$\overset{\swarrow}{=}$}}}
\newcommand{\tw}{\operatorname{tw}}
\newcommand{\base}{\operatorname{base}}
\newcommand{\trans}{\operatorname{trans}}
\newcommand{\rest}{|_}
\newcommand{\bbar}{\overline}
\newcommand{\UML}{\operatorname{\UU\MM\LL}}
\newcommand{\EL}{\mathcal{EL}}
\newcommand{\tsum}{\sideset{}{'}\sum}
\newcommand{\tsh}[1]{\left\{\kern-.9ex\left\{#1\right\}\kern-.9ex\right\}}
\newcommand{\Tsh}[2]{\tsh{#2}_{#1}}
\newcommand{\qeq}{\mathrel{\approx}}
\newcommand{\Qeq}[1]{\mathrel{\approx_{#1}}}
\newcommand{\qle}{\lesssim}
\newcommand{\Qle}[1]{\mathrel{\lesssim_{#1}}}
\newcommand{\simp}{\operatorname{simp}}
\newcommand{\vsucc}{\operatorname{succ}}
\newcommand{\vpred}{\operatorname{pred}}
\newcommand\fhalf[1]{\overrightarrow {#1}}
\newcommand\bhalf[1]{\overleftarrow {#1}}
\newcommand\sleft{_{\text{left}}}
\newcommand\sright{_{\text{right}}}
\newcommand\sbtop{_{\text{top}}}
\newcommand\sbot{_{\text{bot}}}
\newcommand\sll{_{\mathbf l}}
\newcommand\srr{_{\mathbf r}}
\newcommand\geod{\operatorname{\mathbf g}}
\newcommand\mtorus[1]{\boundary U(#1)}
\newcommand\A{\mathbf A}
\newcommand\Aleft[1]{\A\sleft(#1)}
\newcommand\Aright[1]{\A\sright(#1)}
\newcommand\Atop[1]{\A\sbtop(#1)}
\newcommand\Abot[1]{\A\sbot(#1)}
\newcommand\boundvert{{\boundary_{||}}}
\newcommand\storus[1]{U(#1)}
\newcommand\Momega{\omega_M}
\newcommand\nomega{\omega_\nu}
\newcommand\twist{\operatorname{tw}}
\newcommand\modl{M_\nu}
\newcommand\MT{{\mathbb T}}
\newcommand\Teich{{\mathcal T}}
\renewcommand{\Re}{\operatorname{Re}}
\renewcommand{\Im}{\operatorname{Im}}

\title{Caratheodory convergence of log-Riemann surfaces and
Euler's formula}

\author{Kingshook Biswas}
\address{RKM Vivekananda University, Belur Math, WB-711 202, India}

\author{Ricardo Perez-Marco}
\address{CNRS, LAGA, UMR 7539, Universit\'e
Paris 13, Villetaneuse, France}

\begin{abstract}
We define the notion of log-Riemann surfaces
and Caratheodory convergence of log-Riemann surfaces. We prove a
convergence theorem for uniformizations of simply connected
log-Riemann surfaces converging in the Caratheodory topology. We
obtain as a corollary a purely geometric proof of Euler's formula
$\left(1 + \frac{z}{n}\right)^n \to e^z$.
\end{abstract}

\bigskip

\maketitle

\tableofcontents

\section{Introduction}

\medskip

In this article we introduce the notion of {\it log-Riemann
surfaces}. These are Riemann surfaces informally described as
obtained by cutting and pasting copies of $\C$ isometrically. The
larger class of {\it tube-log Riemann surfaces}, allowing pasting
of flat cylinders $\C/\Z$ as well, was used by the second author
for the construction of holomorphic diffeomorphisms with
irrationally indifferent fixed points with special dynamical
properties, see \cite{perezmens}, \cite{perezminvent},
\cite{perezmsmooth}. They were also used by the first author to
construct further examples in \cite{biswas1}, \cite{biswas2}. Our
aim in this article is to initiate a general study of log-Riemann
surfaces.

\medskip

A log-Riemann surface $\mathcal S$ comes equipped with a local
diffeomorphism $\pi : {\mathcal S} \to \C$ satisfying certain
conditions (defined precisely in the next section). The map $\pi$
serves as a coordinate at all points of $\mathcal S$ and allows
one to write formulae for the uniformizations of simply connected
log-Riemann surfaces. The map $\pi$ also induces a flat metric on
$\mathcal S$; the points $p$ added in the completion are
"ramification points" of $\pi$, i.e. $\pi$ restricted to a small
punctured disc neighbourhood of $p$ is a covering of a punctured
disc in $\C$. The order of a ramification point is defined to be
the degree $1 \leq n \leq \infty$ of the covering; the finite
order ramification points can be added to $\SSS$ to give a Riemann
surface ${\mathcal S}^\times$.

\medskip

For a pointed log-Riemann surface $({\mathcal S}^\times, z_0)$ such that
${\mathcal S}^\times$ is simply connected, there is a unique $0 < R \leq
+\infty$, called the {\it conformal radius} of $({\mathcal S}, z_0)$, and
a unique normalized uniformization $F = F_{({\mathcal S}, z_0)} : {\mathcal S}^\times \to
\mathbb{D}_R$ such that $F(z_0) = 0, F'(z_0) =
1$ (with derivative computed in the chart $\pi$).
The flat metric also allows one to define a version of
Gromov-Hausdorff convergence of pointed log-Riemann surfaces,
which we call Caratheodory convergence by analogy with convergence
of simply connected domains in $\C$. We then prove an analogue
of Caratheodory Kernel Convergence Theorem for log-Riemann
surfaces:

\medskip

\begin{theorem} \label{carconvthm} Let $({\mathcal S}_n , z_n) \to ({\mathcal S}, z_0)$ be a
Caratheodory convergent sequence of log-Riemann surfaces such that
the finite completions ${\mathcal S}_n^\times, {\mathcal S}^\times$ are simply
connected. Let the conformal radii and normalized uniformizations
of $({\SSS}_n, z_n), (\SSS, z_0)$ be $R_n, R$ and $F_n, F$
respectively. If $ \limsup_{n\to +\infty} R_n \leq R \ ,$ then
$\lim_{n\to +\infty} R_n = R$ and the uniformizations $F_n$
converge uniformly on compacts of $\SSS$ to the uniformization $F$.
\end{theorem}

\medskip

The convergence of $F_n$ to $F$ above holds in the following sense:
every compact $K \subset \SSS$ containing $z_0$ embeds isometrically
into $\SSS_n$ for $n$ large enough via an isometry $\iota_n$ which is a
translation in the charts $\pi, \pi_n$, and the maps $F_n \circ \iota_n$ converge
uniformly to $F$.

\medskip

We also have convergence of the inverse mappings assuming the
projections of the basepoints converge:

\medskip

\begin{theorem} \label{invconv} Suppose $({\SSS}_n, z_n)$ converge to $({\SSS}, z_0)$. Let $R_n, R$ be the
corresponding conformal radii and $G_n : \mathbb{D}_{R_n} \to {\SSS}_n^\times, G : \mathbb{D}_R \to
{\SSS}^\times$ be normalized biholomorphisms such that $G_n(0) =
z_n, G(0) = z_0$ and $(\pi_n \circ G_n)'(0) = 1, (\pi \circ G)'(0)
= 1$. If $\limsup R_n \leq R$ and $\pi_n(z_n) \to \pi(z_0)$ then
$\lim R_n = R$ and $\pi_n \circ G_n \to \pi \circ G$ uniformly on
compacts of $\mathbb{D}_R$.
\end{theorem}

\medskip

Observing that for the log-Riemann surfaces $\SSS_n, \SSS$
corresponding to the $n$th root and logarithm functions
respectively, with base point $1$, the normalized uniformizations are
given in charts by the functions $\pi_n \circ G_n(z) = (1 + \frac{z}{n})^n$
and $\pi \circ G(z) = e^z$ respectively, we obtain as a direct corollary of the
above a purely geometric proof of Euler's formula:

\medskip

\begin{cor} \label{euler} We have
$$
\left( 1 + \frac{z}{n} \right)^n \to e^z
$$
uniformly on compacts.
\end{cor}

\medskip

In a forthcoming article \cite{bipmlogrs1} we use the above
convergence theorems to obtain uniformization formulae for simply
connected log-Riemann surfaces with finitely many ramification
points.

\medskip

\subsection{Acknowledgements} The authors are very grateful to
J.Mu\~{n}oz D\'{i}az, R. Alonso Blanco, J. Lombardero, S.
Jim'enez Verdugo, A. \'{A}lvarez V\'{a}zquez, D. Barsky
and F. Marcellan Espa\~{n}ol for numerous helpful discussions. The
first author was partly supported by Department of Science and Technology research project
grant DyNo. 100/IFD/8347/2008-2009, and the second author by CNRS, UMR 7539.

\bigskip

\section{Log-Riemann surfaces}

\medskip

\subsection{Definitions}

\medskip

\begin{definition} A log-Riemann surface is a
Riemann surface $\SSS$ equipped with a local diffeomorphism $\pi
: \SSS \to \C$, called the projection mapping, such that the
following holds: let $\overline{\SSS} = \SSS \coprod \RR$ be the
completion of $\SSS$ with respect to the path-metric on $\SSS$
induced by the flat metric $|d\pi|$ (the pull-back of the
Euclidean metric under $\pi$). Then $\RR$ is discrete.

The set $\RR$ is called the ramification locus of $\SSS$ and
the points of $\RR$ are called ramification points.
\end{definition}

\medskip

Since $\pi$ is a local isometry with respect to the metrics on
$\SSS$ and $\C$ it extends uniquely to $\overline{\SSS}$.

\medskip

\begin{prop} Let $p \in \RR$.
For $r > 0$ small enough, $\pi(B(p,r) - \{p\}) = B(\pi(p), r) - \{\pi(p)\}$
and $$\pi : B(p, r) - \{p\} \to
B(\pi(p), r) - \{\pi(p)\}$$ is a covering.
\end{prop}

\medskip

\noindent{\bf Proof:} Choose $r>0$ small enough so that there are no points of $\RR$
 in $B(p,r)$ apart from $p$. Pick a $z \in B(p,r) - \{p\}$;
then $\pi(z) \in B(\pi(p),r)$, and the local inverse
of $\pi$ at $\pi(z)$ such that $\pi^{-1}(\pi(z)) = z$ can be
analytically continued to all points of $B(\pi(p),r) - \{\pi(p)\}$,
since the only possible obstruction to the continuation is
encountering points in $\RR - \{p\}$ above, but by the choice of $r$ this
is not possible. Therefore $\pi$ maps $B(p,r) - \{p\}$ onto
$B(\pi(p),r) - \{\pi(p)\}$. Moreover the local inverse is
single-valued on $B(\pi(z), t)$ where $t = d(z,p)$, so $\pi$ maps
$B(z,t)$ isometrically to $B(\pi(z),t)$, so taking $z_n \in
B(z,t)$ converging to $p$, we have
$$
|\pi(z) - \pi(p)| = \lim |\pi(z) - \pi(z_n)| = \lim d(z, z_n) =
d(z,p) > 0
$$
so $\pi(z) \neq \pi(p)$.


The proof that $\pi: B(p,r)-\{p\}\to \pi (B(p,r)-\{p\})$
is a covering is similar. For each point
$z_0 \in \pi (B(p,r)-\{p\})=B(\pi(p),r)-\{\pi(p)\}$,
take a small disc $B(z_0,\rho) \subset B(\pi(p), r)-\{
\pi(p)\}$ and let $U$ be a connected component of
$\pi^{-1}(B(z_0, \rho))$; we can pick a $z_1 \in U$ and as before
continue without obstruction the local inverse of $\pi$ at $\pi(z_1)$ satisfying
$\pi^{-1}(\pi(z_1)) = z_1$ to all of the disk $B(z_0,\rho)$; the
continuation of $\pi^{-1}$ to it
is single-valued, so $\pi_{|U} : U \to B(z_0,\rho)$ is a diffeomorphism. $\diamond$

\medskip

\begin{definition} The order $1 \leq
n \leq \infty$ of a ramification point $p \in \RR$ is defined to
be the degree of the covering $\pi : B(p, r) - \{p\} \to B(\pi(p),
r) - \{\pi(p)\}$ of the punctured disk $B(\pi(p), r) - \{\pi(p)\}$
(where $r$ is taken small enough as above). The ramification point $p$ is
called finite if $n < \infty$ and infinite if $n = \infty$.
\end{definition}

\medskip

\begin{definition} The finite completion ${\SSS}^\times \subset
\overline{\SSS}$ of a
log-Riemann surface $\SSS$ is defined to be the union of $\SSS$ and all
finite ramification points of $\SSS$.
\end{definition}

\medskip

Since any finite sheeted holomorphic covering of a punctured disc is
equivalent to the covering $\mathbb{D}^* \to \mathbb{D}^*, z \mapsto z^n$, a punctured
neighbourhood of a finite ramification point is
biholomorphic to a punctured disk, so the finite ramification points may be added to
$\SSS$ to give a Riemann surface structure on ${\SSS}^\times$ compatible with that
of $\SSS$, such that the map $\pi : {\SSS}^\times \to \C$ is holomorphic with critical points
at the finite ramification points.

\medskip

When the finite completion is simply connected, ${\SSS}^\times$ is
biholomorphic to $\mathbb{D}_R = \{ |z| < R \}$ for
some $0 < R \leq \infty$. Fixing a basepoint $z_0 \in \SSS$, we
can choose $R$ and a uniformization $G : \mathbb{D}_R
\to \SSS$ satisfying $G(0) = z_0, (\pi \circ G)'(0) = 1$, and it is
easy to see that $R$ and $G$ are then uniquely determined.

\medskip

\begin{definition} For $(\SSS, z_0)$ a pointed
log-Riemann surface such that ${\SSS}^\times$ is simply connected,
the conformal radius $R$ of $(\SSS, z_0)$ is the unique $0 < R
\leq \infty$ as determined above.
\end{definition}

\medskip

\subsection{Examples}

\medskip

We briefly describe some examples of log-Riemann surfaces:

\medskip

\noindent {\bf 0. Log-Riemann surfaces associated to the complex plane:} Let
$\SSS = \C $ and $\pi : \SSS \to \C$ be any automorphism of the plane $\pi(z) = az+b$, 
$a\in \C^*$, $b\in \C$. The associated metric space is the euclidean plane with the euclidean metric
scaled by $|a|$. The ramification locus $\RR$ is empty. It is easy to see that these are 
the only (connected) log-Riemann surfaces with $\RR=\emptyset$.

\medskip

\noindent {\bf 1. The log-Riemann surface of the nth root:} Let
$\SSS = \C - \{0\}$ and $\pi : \SSS \to \C, \pi(z) = z^n$.
It is easy to see that $(\SSS, \pi)$ is a log-Riemann surface isometric to the
union of $n$ slit planes $\C - [-\infty,0]$
pasted isometrically along the $2n$ ''sides'' of the slits, with a single ramification
point $p$ of order $n$, $\overline{\SSS} = {\SSS}^\times = \SSS \cup
\{p\}$, such that $\pi(p) = 0$.

\medskip

\noindent {\bf 2. Polynomial log-Riemann surfaces:} Generalizing
the previous example, let $\pi$ be a polynomial with $C$ the set
of critical points. Then $\SSS = \C - C$ is a log-Riemann
surface with $\#C$-many ramification points whose orders add up to
the degree of $\pi$.

\medskip

\noindent {\bf 3. The log-Riemann surface of the logarithm:} Let
$\SSS = \C - \{0\}$ and $\pi : \SSS \to \C, \pi(z) = e^z$.
It is easy to see that $(\SSS, \pi)$ is a log-Riemann surface isometric to the
union of infinitely many slit planes $\C - [-\infty,0]$
pasted isometrically along the slits, with a single ramification
point $p$, of infinite order, $\overline{\SSS} = \SSS \cup
\{p\}$, such that $\pi(p) = 0$.

\medskip

\noindent {\bf 4. The universal covering of a log-Riemann
surface:} Let $\pi : \tilde{\SSS} \to \SSS$ be the universal
covering of a log-Riemann surface $\SSS$ with projection mapping
$\pi_{\SSS} : \SSS \to \C$. The map $\pi_{\tilde{\SSS}} :=
\pi_{\SSS} \circ \pi$ is a local holomorphic diffeomorphism. With respect
to the induced path metric on $\tilde{\SSS}$ the map $\pi$ is a $1$-Lipschitz
local isometry and hence extends to a map between the completions
$\pi : \tilde{\SSS}^* \to \SSS^*$. Moreover since $\RR$ is discrete,
$\tilde{\SSS}^* - \tilde{\SSS} = \pi^{-1}(\RR)$ is discrete, and
hence endows $\tilde{\SSS}$ with a
log-Riemann surface structure.

\bigskip

\section {Caratheodory Theorem for log-Riemann surfaces}

\medskip

\subsection{Caratheodory convergence of log-Riemann surfaces}

\medskip

Recall that log-Riemann surfaces are endowed with the flat metric
$|d\pi|$.

\medskip

\begin{definition} A pointed sequence of log-Riemann surfaces
$({\SSS}_n ,z_n )$ converges to a pointed log Riemann surface $({\SSS}, z_0)$
if for any compact $K \subset \SSS$ containing $z_0$
there exists $N=N(K)\geq 1$ such that for $n\geq
N$ there is an isometric embedding $\iota$ of $K$ into ${\SSS}_n$
mapping $z_0$ into $z_n$ such that $\iota$ is a translation in the
charts $\pi, \pi_n$ (the translation that maps $\pi(z_0)$ to $\pi(z_n)$).
\end{definition}

\medskip

\noindent {\bf Example.} Let
$(\SSS , z_0)$ be a pointed log-Riemann surface 
with projection mapping $\pi : \SSS \to \C$, $\pi (z_0)=0$. 
For $n\geq 1$ we consider the new log-Riemann surface $\SSS_n$ 
with projection mapping $\pi_n = n \cdot \pi$ and $z_n=z_0$.
When $n\to +\infty$, we have that $(\SSS_n, z_n)$ has the planar log-Riemann surface 
as Caratheodory limit.

\medskip

\medskip

The Caratheodory limit of a sequence $({\SSS}_n, z_n)$ if it exists is unique
up to isometry:

\medskip

\begin{prop} Let $({\SSS}, z_0)$ and $({\SSS}', z_0')$ be two
pointed log-Riemann surfaces both of which are Caratheodory limits
of a sequence of pointed log-Riemann surfaces $({\SSS}_n, z_n)$.
Then there is an isometry $T : \SSS \to
{\SSS}'$ taking $z_0$ to $z_0'$ whose
expression in log-charts is the translation mapping $\pi(z_0)$
to $\pi'(z_0')$.
\end{prop}

\medskip

\noindent{\bf Proof:} Consider the germ of holomorphic diffeomorphism $T$
from $\SSS$ to ${\SSS}'$ mapping $z_0$ to $z_0'$ whose
expression in log-charts is the translation maping $\pi(z_0)$
to $\pi'(z_0')$.

\medskip

Let $\gamma : [a,b] \to \SSS$ be a
curve in $\SSS$ starting from $z_0$ along
which $T$ can be continued. For $n$ large enough, if $\iota$
and $\iota'$ denote the isometric embeddings of $\gamma$ and
$T(\gamma)$ respectively into ${\SSS}_n$, on $\gamma$ we must
have
$$
\iota = \iota' \circ T
$$

\medskip

\begin{lemma} The germ $T$ can be continued analytically along all
paths in $\SSS$.
\end{lemma}

\medskip

\noindent{\bf Proof:} Suppose there is a path $\gamma : [a,b] \to \SSS,
\, \gamma(a) = z_0, \gamma(b) = z \in \SSS$ such that $T$ can be continued
analytically along $\gamma([a,b))$ but not up to $\gamma(b) = z$.
Since $T$ is a local isometry, the limit $z' = \lim_{x \to b} T(\gamma(x))$
exists and must be a ramification point of order $n > 1$ since $T$ cannot be
continued to $\gamma(b)$.

\medskip

Take $\delta > 0$ small enough so that $B(z',\delta)$ contains no
other ramification points and so that $\overline{B(z,\delta)}
\subset \SSS$. Let $b_1 \in [a,b)$ be such that $T(\gamma(b_1))
\in B(z',\delta)$. Let $\alpha : [c,d] \to \SSS$ be a circular
loop in $\SSS$ that winds once around $z$, starting from
$\gamma(b_1)$, $\alpha(c) = \gamma(b_1)$. We note that $T$ can be
continued along $\alpha$, but $T(\alpha(c))$ is not equal to
$T(\alpha(d))$.

\medskip

For $n$ large enough the compacts $\gamma([a,b])
\cup \alpha([c,d]) \cup \overline{B(z,\delta)} \subset \SSS$ and
$T(\gamma([a,b_1])) \cup T(\alpha([c,d])) \subset {\SSS}'$ both
embed isometrically into ${\SSS}_n$, let $\iota, \iota'$ be
the respective embeddings. Now, the ball $\iota(B(z,\delta))$ is
completely contained in ${\SSS}_n$, so for the curve $\alpha$ we
have
$$
\iota(\alpha(c)) = \iota(\alpha(d)) = \iota(\gamma(b_1)) \ \hbox{
where } \alpha(c) = \alpha(d) = \gamma(b_1),
$$
which implies
$$
\iota'(T(\alpha(c))) = \iota'(T(\alpha(d)))
$$
and hence, since $\iota'$ is an isometry,
$$
T(\alpha(c)) = T(\alpha(d)),
$$
a contradiction. $\diamond$

\medskip

\begin{lemma} The continuation of $T$ to all of $\SSS$ is
single-valued.
\end{lemma}

\medskip

\noindent{\bf Proof:} Let $\gamma : [a,b] \to \SSS, \, \gamma(a) =
\gamma(b) = z_0$,
 be a closed path in $\SSS$. Consider the curve $T(\gamma) \in {\SSS}'$ given
by continuing $T$ along $\gamma$. Take $n$ large enough so that
the compacts $\gamma([a,b]) \subset \SSS$ and $T(\gamma([a,b]))
\subset {\SSS}'$ both embed isometrically into ${\SSS}_n$ via
isometries $\iota$ and $\iota'$ respectively. As before, along
$\gamma$ we have
$$
\iota = \iota' \circ T
$$
Since
$$
\iota(\gamma(a)) = \iota(\gamma(b)) = \iota(z_0) = z_n,
$$
it follows that
$$
\iota'(T(\gamma(a))) = \iota'(T(\gamma(b)))
$$
and hence
$$
T(\gamma(a))) = T(\gamma(b))),
$$
$\diamond$

\medskip

It follows from the above lemmas that we obtain a globally defined
map $T : \SSS \to {\SSS}'$. Applying the same arguments to the
germ $S = T^{-1}$ given by the inverse of the initial germ $T$
gives a map $S : {\SSS}' \to \SSS$, and it is straightforward to
check that $T$ and $S$ define global mutual inverses. The
conclusion of the Proposition follows. $\diamond$

\bigskip

\subsection{Convergence of uniformizations}

\medskip

\noindent{\bf Proof of Theorem \ref{carconvthm} :} Let $C \subset
\mathbb{D}_R$ be the image under $F$ of the finite ramification
points of $\SSS$. Then $C$ is discrete and the complement
$\mathbb{D}_R - C$ contains a ball around the origin. For any
compact $K \subset \SSS$ containing $z_0$, for $n$ large enough
let $\iota_n : K \to {\SSS}_n$ be the corresponding embeddings.
The maps $g_n = F_n \circ \iota_n \circ F^{-1}$ are well-defined
and univalent on any compact in $\mathbb{D}_R - C$ for $n$ large enough,
satisfying $g_n(0) = 0, g_n'(0) = 1$. It follows from classical
univalent function theory that they form a normal family on any
simply connected subdomain of $\mathbb{D}_R - C$ containing the
origin. Since normality is a local property and every point in
$\mathbb{D}_R- C$ has a neighbourhood contained in such a simply
connected subdomain, it follows that the family $(g_n)$ is normal
on $\mathbb{D}_R - C$. Let $g$ be a normal limit on $\mathbb{D}_R
- C$ of the $g_n$'s; then $g$ is univalent and $g(0) = 0$. Hence
for any Jordan curve $\gamma \subset \mathbb{D}_R - C$ if $0$ is
in the bounded component of $\C - \gamma$ then $0$ is in the
bounded component of $\C - g(\gamma)$. It follows that $g$ is
bounded on the bounded component of $\C - \gamma$, and hence the
singularities of $g$ at the points of $C$ are removable.

\medskip

It is easy to see that the extension of $g$ to $\mathbb{D}_R$ is
also univalent; since $g(0) = 0, g'(0) = 1$, if $R = \infty$ then
$g$ must be the identity. If $R < \infty$, since $\limsup R_n \leq
R$, $g$ takes values in $\overline{\mathbb{D}_R}$, so by the
Schwarz lemma $g$ is the identity. So the only limit point of the
sequence $g_n$ is the identity, and the sequence converges to the
identity. Hence $F_n$ converges to $F$ uniformly on compacts, in
the sense that $F_n \circ \iota_n = g_n \circ F$ converges on
compacts of $\SSS$ to $F$. If $\liminf R_n < R$, then along a
subsequence the maps $g_n$ take values in a disc compactly
contained in $\mathbb{D}_R$, contradicting the convergence of
$g_n$ to the identity; hence $\liminf R_n \geq R$, so $\lim R_n =
R$. $\diamond$

\medskip

\noindent{\bf Proof of Theorem \ref{invconv}:} In the notation as
above, $G_n = F_n^{-1}, G = F^{-1}$, and we have seen that the
univalent maps $g_n = F_n \circ \iota_n \circ F^{-1} \to id$
uniformly on compacts in $\mathbb{D}_R - C$, $\lim R_n = R$. It
follows that $g_n^{-1} \to id$ and so $(F_n \circ \iota_n)^{-1} =
G \circ g_n^{-1} \to G$ uniformly on compacts in $\mathbb{D}_R -
C$. Also for any compact in $\mathbb{D}_R - C$ for $n$ large we can
write
$$
\pi_n \circ G_n = (\pi_n \circ \iota_n) \circ (F_n \circ \iota_n)^{-1}
$$
By the hypothesis on the isometric embeddings,
$\pi_n \circ \iota_n = T_n \circ \pi$ where $T_n$ is the
translation sending $\pi(z_0)$ to $\pi_n(z_n)$; since $\pi_n(z_n)
\to \pi(z_0)$, $T_n \to id$, so $\pi_n \circ \iota_n \to \pi$
uniformly on compacts in $\SSS$, hence $\pi_n \circ G_n \to \pi
\circ G$ uniformly on compacts in $\mathbb{D}_R - C$, and hence by
the maximum principle also on $\mathbb{D}_R$. $\diamond$

\medskip

\noindent{\bf Proof of Corollary \ref{euler}:} Let ${\SSS}_n$ be the
log-Riemann surface of the $n$th root with projection $\pi_n : {\SSS}_n \to \C$,
and take basepoints $z_n$ such that $\pi_n(z_n) = 1$. Let $\SSS$
be the log-Rieman surface of the logarithm with projection $\pi :
\SSS \to \C$ and basepoint $z_0$ such that $\pi(z_0) = 1$.
Then any compact $K$ in $\SSS$ can only intersect finitely many
''sheets'' of $\SSS$ and hence embeds isometrically in ${\SSS}_n$
for $n$ large enough, so $({\SSS}_n, z_n) \to (\SSS, z_0)$. The
normalized uniformizations $G_n, G$ of ${\SSS}_n, \SSS$ are given
in charts by $\pi_n \circ G_n(z) = (1 + \frac{z}{n})^n$ and $\pi
\circ G(z) = e^z$ respectively so the previous Theorem gives the
conclusion of the Corollary. $\diamond$

%
%
%

\bigskip

\bibliography{carconv}

\begin{thebibliography}{BPM10}

\bibitem[Bis05]{biswas1}
K.~Biswas.
\newblock Smooth combs inside hedgehogs.
\newblock {\em DISCRETE AND CONTINUOUS DYNAMICAL SYSTEMS, vol. 12, 5}, pages
  853--880, 2005.

\bibitem[Bis08]{biswas2}
K.~Biswas.
\newblock Hedgehogs of hausdorff dimension one.
\newblock {\em Ergodic Theory and Dynamical Systems, vol. 28, 6}, pages
  1713--1727, 2008.

\bibitem[BPM10]{bipmlogrs1}
K.~Biswas and R.~Perez-Marco.
\newblock Uniformization of simply connected finite type log-riemann surfaces.
\newblock {\em preprint}, 2010.

\bibitem[PM93]{perezmens}
R.~Perez-Marco.
\newblock Sur les dynamiques holomorphes non-linearisables et une conjecture de
  v.i. arnold.
\newblock {\em Annales Scientifiques de l'E.N.S 26}, pages 565--644, 1993.

\bibitem[PM95]{perezminvent}
R.~Perez-Marco.
\newblock Uncountable number of symmetries for non-linearisable holomorphic
  dynamics.
\newblock {\em Inventiones Mathematicae 119}, pages 67--127, 1995.

\bibitem[PM00]{perezmsmooth}
R.~Perez-Marco.
\newblock Siegel disks with smooth boundary.
\newblock {\em preprint}, 2000.

\end{thebibliography}
\bibliographystyle{alpha}

\end{document}